\def\titlep{A tensor product of representations of Cuntz algebras}
\documentclass[8pt]{article}
\usepackage{graphicx,ifthen}
\usepackage{amssymb}

\font\germ=eufm10 at12pt

\def\goth#1{\hbox{\germ#1}}

\newcommand{\qed}{\hbox{\rule[-2pt]{3pt}{6pt}}}
\newcommand{\qedh}{\hfill\qed \\}

\newcommand{\vv}{\vspace{.3in}}

\newcommand{\vep}{\varepsilon}

\newtheorem{Thm}{Theorem}[section]

\newtheorem{defi}[Thm]{Definition}
\newtheorem{lem}[Thm]{Lemma}

\newtheorem{prop}[Thm]{Proposition}

\newtheorem{cor}[Thm]{Corollary}

\newcommand{\ww}{\vv\noindent}

\def\cal#1{\mathcal #1}
\def\con{{\cal O}_{N}}

\def\edot{=1,\ldots,N}
\def\pr{{\it Proof.}\quad}

\def\nset#1{\{1,\ldots,N\}^{#1}}

\def\co#1{{\cal O}_{#1}}
\def\disp#1{{\displaystyle #1}}
\def\brl{branching law}
\def\bfsnl{{\rm BFS}_{N}(\Lambda)}

\addtocounter{footnote}{1}
\def\sftt#1{
\setcounter{equation}{0}
\addtocounter{footnote}{1}
\section{#1}
}

\def\ssft#1{\subsection{#1}}

\begin{document}
%
%
\def\autherp{Katsunori Kawamura}
\def\emailp{E-mail: kawamura@kurims.kyoto-u.ac.jp.}
\def\addressp{College of Science and Engineering Ritsumeikan University,\\
1-1-1 Noji Higashi, Kusatsu, Shiga 525-8577,Japan
}

\def\infw{\Lambda^{\frac{\infty}{2}}V}
\def\zhalfs{{\bf Z}+\frac{1}{2}}
\def\ems{\emptyset}
\def\pmvac{|{\rm vac}\!\!>\!\! _{\pm}}
\def\vac{|{\rm vac}\rangle _{+}}
\def\dvac{|{\rm vac}\rangle _{-}}
\def\ovac{|0\rangle}
\def\tovac{|\tilde{0}\rangle}
\def\expt#1{\langle #1\rangle}
\def\zph{{\bf Z}_{+/2}}
\def\zmh{{\bf Z}_{-/2}}
\def\brl{branching law}
\def\bfsnl{{\rm BFS}_{N}(\Lambda)}
\def\scm#1{S({\bf C}^{N})^{\otimes #1}}
\def\mqb{\{(M_{i},q_{i},B_{i})\}_{i=1}^{N}}
\def\zhalf{\mbox{${\bf Z}+\frac{1}{2}$}}
\def\zmha{\mbox{${\bf Z}_{\leq 0}-\frac{1}{2}$}}
\newcommand{\mline}{\noindent
\thicklines
\setlength{\unitlength}{.1mm}
\begin{picture}(1000,5)
\put(0,0){\line(1,0){1250}}
\end{picture}
\par
 }
\def\ptimes{\otimes_{\varphi}}
\newcommand{\bolda}{\hbox{\boldmath$\mathit{a}$}}
\newcommand{\boldb}{\hbox{\boldmath$\mathit{b}$}}
\newcommand{\boldx}{\hbox{\boldmath$\mathit{x}$}}
\def\pca{pseudo-Cuntz algebra}
%
%
%
\setcounter{section}{0}
\setcounter{footnote}{0}
\setcounter{page}{1}
\pagestyle{plain}

%
%

\noindent
{\Large \titlep}

\ww
\autherp

\noindent
{\it \addressp}\\
{\it \emailp}\\
\quad \\


%
%

\noindent
{\bf Abstract.}
We introduce a nonsymmetric, associative tensor product 
among representations of Cuntz algebras by using embeddings.
We show the decomposition formulae of tensor products
for permutative representations explicitly 
We apply decomposition formulae to determine
properties of endomorphisms.
\\


%
%
\sftt{Introduction}
\label{section:first}
%
%
\ssft{Motivation, definition and basic facts}
\label{subsection:firstone}
In algebraic quantum field theory,
endomorphisms of operator algebras play important role \cite{DHR}. 
We have shown in \cite{PE01,PE04, PE02},
the branching laws for endomorphisms of the Cuntz algebra $\con$
with respect to permutative representations by \cite{BJ,DaPi2,DaPi3},
and shown properties and classifications of endomorphisms
by branching laws.

When I gave a talk on this in Kyoto, September, 2004, 
I. Ojima asked me a specific question of what is 
the tensor product of representations of $\con$.
This is indeed impossible to answer in the usual sense
because nobody knows the ``canonical" (or suitable)
embedding of $\con$ into $\con\otimes \con$. 
By \cite{Rordam}, it is known that $\co{2}\cong \co{2}\otimes\co{2}$.
However it seems that the property of such map is not easy to
study because we can not find concrete formulae of 
this isomorphism with respect to canonical generators of $\co{2}$.
Hence we give up to use such isomorphism in this study and
we consider other possibility. 
For this purpose, we begin to review the
definition of tensor product in group theory.

By the diagonal embedding $\varphi$ of a group $G$ into $G\times G$, 
we obtain the representation $\pi_{1}\ptimes \pi_{2}$ of $G$ by
\[\pi_{1}\ptimes \pi_{2}\equiv (\pi_{1}\otimes \pi_{2})\circ \varphi\]
from the tensor representation $\pi_{1}\otimes \pi_{2}$ of $G\times G$ 
for representations $\pi_{1}$ and $\pi_{2}$ of $G$.
The representation $\pi_{1}\ptimes \pi_{2}$
is usually called the tensor product of $\pi_{1}$ and $\pi_{2}$ simply.
The associativity and the distribution law with respect to
the direct sum for $\ptimes$ is assured by the property of $\varphi$. 
Because $\pi_{1}\ptimes \pi_{2}$ and 
$\pi_{2}\ptimes \pi_{1}$ are equivalent, this tensor product 
is called {\it symmetric} (or commutative).
If there is other embedding $\varphi^{'}$ of $G$ into $G\times G$,
we can define another composition of two representations of $G$ by $\varphi^{'}$
and its property also depends on the choice of $\varphi^{'}$.
In this sense, there may exist other tensor like structure
among representations of $G$.
According to this idea, we define a kind of tensor product 
of representations of Cuntz algebras as follows.

%
For $2\leq N,M<\infty$,
let $t_{1},\ldots,t_{N},r_{1},\ldots,r_{M}$ and $s_{1},\ldots,s_{NM}$
be canonical generators of $\con$, $\co{M}$ and $\co{NM}$, respectively.
Define the unital 
$*$-embedding $\varphi_{N,M}$ of $\co{NM}$ into $\con\otimes \co{M}$ by
%
%
\begin{equation}
\label{eqn:standard}
\varphi_{N,M}(s_{M(i-1)+j})\equiv t_{i}\otimes r_{j}\quad(i\edot,\,j=1,\ldots,M).
\end{equation}
Then the following diagram is commutative for each $N,M,L\geq 2$:

\noindent
%
%
%
\thicklines
\setlength{\unitlength}{.1mm}
\begin{picture}(1000,300)(-50,120)
\put(0,250){$\co{NML}$}
\put(150,280){\vector(3,1){200}}
\put(150,350){$\varphi_{N,ML}$}
\put(150,240){\vector(3,-1){200}}
\put(150,160){$\varphi_{NM,L}$}
\put(400,350){$\co{N}\otimes \co{ML}$}
\put(400,150){$\co{NM}\otimes \co{L}$}
\put(650,350){\vector(3,-1){200}}
\put(700,350){$id_{N}\otimes \varphi_{M,L}$}
\put(650,170){\vector(3,1){200}}
\put(700,160){$\varphi_{N,M}\otimes id_{L}$}
\put(900,250){$\co{N}\otimes \co{M}\otimes \co{L}$}
\end{picture}

\noindent
In other words, the following holds on $\co{NML}$:
%
%
\begin{equation}
\label{eqn:coassociativity}
(id_{N}\otimes \varphi_{M,L})\circ \varphi_{N,ML}
=(\varphi_{N,M}\otimes id_{L})\circ \varphi_{NM,L}.
\end{equation}
We call this property the {\it weak co-associativity} 
of the family $\varphi\equiv \{\varphi_{N,M}:N,M\geq 2\}$.
Define ${\rm Rep}\con$ the class of all unital $*$-representations of $\con$.
For the family $\varphi$,
we introduce the following operation $\ptimes$.
%
%
\begin{defi}
\label{defi:fundamental}
For $(\pi_{1},\pi_{2})\in {\rm Rep}\con\times {\rm Rep}\co{M}$,
the $\varphi$-tensor product 
$\pi_{1}\otimes_{\varphi} \pi_{2}\in{\rm Rep}\co{NM}$ 
of $\pi_{1}$ and $\pi_{2}$ is  defined by 
%
%
\begin{equation}
\label{eqn:definition}
\pi_{1}\ptimes\pi_{2}\equiv(\pi_{1}\otimes \pi_{2})\circ \varphi_{N,M}.
\end{equation}
\end{defi}

\noindent
For $\pi_{1},\pi_{2}\in{\rm Rep}\con$,
we denote $\pi_{1}\sim \pi_{2}$ if $\pi_{1}$ and $\pi_{2}$ are equivalent.
%
%
\begin{lem}
\label{lem:fundamental}
For $\pi_{1},\pi_{1}^{'}\in {\rm Rep}\con$,
$\pi_{2},\pi_{2}^{'}\in {\rm Rep}\co{M}$ and $\pi_{3}\in {\rm Rep}\co{L}$,
the following holds.
\begin{enumerate}
\item
If $\pi_{1}\sim \pi_{2}^{'}$ and $\pi_{2}\sim \pi_{2}^{'}$,
then $\pi_{1}\otimes_{\varphi} \pi_{2}\sim
\pi_{1}^{'}\otimes_{\varphi} \pi_{2}^{'}$.
\item
$\pi_{1}\otimes_{\varphi} (\pi_{2}\oplus \pi_{2}^{'})=
\pi_{1}\otimes_{\varphi} \pi_{2}\,\oplus\, \pi_{1}\otimes_{\varphi} \pi_{2}^{'}$.
\item
$\pi_{1}\otimes_{\varphi} (\pi_{2}\otimes_{\varphi} \pi_{3})
=(\pi_{1}\otimes_{\varphi} \pi_{2})\otimes_{\varphi} \pi_{3}$.
\end{enumerate}
\end{lem}

\noindent
Lemma \ref{lem:fundamental} (i) and (ii) are easily verified.
The statement (iii) is derived from the weak co-associativity of $\varphi$
in (\ref{eqn:coassociativity}).
By Lemma \ref{lem:fundamental} (iii),
we can use a notation $\pi_{1}\otimes_{\varphi}\cdots
\otimes_{\varphi}\pi_{n}$ for $\pi_{1},\ldots,\pi_{n}$.

We show a relation between the permutation of the order of the tensor product
and automorphisms as follows.
For $N,M\geq 2$,
define the bijection $\phi_{N,M}$ from 
$X\equiv \{1,\ldots,N\}\times\{1,\ldots,M\}$ onto $\{1,\ldots,NM\}$ by
%
%
\begin{equation}
\label{eqn:bijection}
\phi_{N,M}(a,b)\equiv M(a-1)+b\quad((a,b)\in X).
\end{equation}
%
%
\begin{Thm}
\label{Thm:symmetry}
For $N_{1},\ldots,N_{n}\geq 2$ with $n\geq 2$,
let $M\equiv N_{1}\cdots N_{n}$ and
$Y_{i}\equiv \{1,\ldots,N_{i}\}$ for $i=1,\ldots,n$.
For a permutation $\sigma\in {\goth S}_{n}$, define
\[\tau_{\sigma,i}\equiv \phi_{N_{\sigma(1)}\cdots N_{\sigma(i-1)},
N_{\sigma(i)}}\quad(i=2,\ldots,n.)\]
For $(i_{1},\ldots,i_{n})\in {\cal Y}\equiv Y_{1}\times \cdots \times Y_{n}$ 
and $\sigma\in {\goth S}_{n}$, define
$[i_{1},\ldots,i_{n}]_{\sigma}\in\{1,\ldots,M\}$ by
\[[i_{1},\ldots,i_{n}]_{\sigma}
\equiv \{\tau_{\sigma,n}\circ (\tau_{\sigma,n-1}\times id)
\circ\cdots \circ(\tau_{\sigma,2}\times id)\}
(i_{\sigma(1)},\ldots,i_{\sigma(n)}).\]
For $\sigma,\eta\in {\goth S}_{n}$,
define the automorphism $\alpha_{\sigma,\eta}$ of $\co{M}$ by 
\[\alpha_{\sigma,\eta}(s_{[i_{1},\ldots,i_{n}]_{\eta}})\equiv
s_{[i_{1},\ldots,i_{n}]_{\sigma}}\quad
((i_{1},\ldots,i_{n})\in {\cal Y}).\]
Then for any $\pi_{i}\in{\rm Rep}\co{N_{i}}$ with
$i=1,\ldots,n$, the following holds:
\[\pi_{\sigma(1)}\ptimes\cdots\ptimes\pi_{\sigma(n)}\sim
(\pi_{\eta(1)}\ptimes\cdots\ptimes\pi_{\eta(n)})
\circ \alpha_{\eta,\sigma}.\]
\end{Thm}

We summarize results here.
Remark that $\pi_{1},\pi_{2},\pi_{1}\otimes_{\varphi}\pi_{2}$
in Definition \ref{defi:fundamental} are representations of 
different algebras in general.
This differs from the tensor product in group theory.
By Lemma \ref{lem:fundamental} (i),
$\otimes_{\varphi}$ is well-defined on the set of equivalence classes. 
Especially $\pi_{1}\otimes_{\varphi}\pi_{2}$ is unique up to
unitary equivalence.
By Lemma \ref{lem:fundamental} (ii) and (iii),
$\otimes_{\varphi}$ is associative and
distributive with respect to the direct sum.
By Theorem \ref{Thm:symmetry},
$\pi_{1}\otimes_{\varphi}\pi_{2}$ and $\pi_{2}\otimes_{\varphi}\pi_{1}$
are not equivalent in general
because $\alpha_{id,\sigma}$ is outer when $\sigma\ne id$.
This is different from the situation of quantum groups with $R$-matrix
\cite{Drinfeld}.
We show concrete examples such that 
$\pi_{1}\otimes_{\varphi}\pi_{2}\not\sim \pi_{2}\otimes_{\varphi}\pi_{1}$
in $\S$ \ref{subsection:fourthone}.
It seems that the definition of $\otimes_{\varphi}$ is artificial,
but $\otimes_{\varphi}$ satisfies ingredients of tensor product
except its commutativity. 
Furthermore, the decomposition formulae associated with $\ptimes$ 
is explicitly computed on permutative representations in the next subsection.

%
%
\ssft{Tensor product of permutative representations}
\label{subsection:firsttwo}
In this paper, any representation, embedding and endomorphism
are assumed unital and $*$-preserving.
%
%
\begin{defi}
\label{defi:first}
Let $s_{1},\ldots,s_{N}$ be canonical generators of $\con$ and
let $({\cal H},\pi)$ be a representation of $\con$.
\begin{enumerate}
\item
$({\cal H},\pi)$ is a {\it permutative representation} of $\con$
if there is a complete orthonormal basis $\{e_{n}\}_{n\in\Lambda}$
of ${\cal H}$ and a family $f=\{f_{i}\}_{i=1}^{N}$
of maps on $\Lambda$ such that $\pi(s_{i})e_{n}=e_{f_{i}(n)}$
for each $n\in\Lambda$ and $i\edot$.
\item
For $J=(j_{l})_{l=1}^{k}\in\nset{k}$,
$P_{N}(J)$ is a class of $({\cal H},\pi)$ with 
a unit cyclic vector $\Omega\in {\cal H}$
such that $\pi(s_{J})\Omega=\Omega$ 
and $\{\pi(s_{j_{l}}\cdots s_{j_{k}})\Omega\}_{l=1}^{k}$
is an orthonormal family in ${\cal H}$
where $s_{J}\equiv s_{j_{1}}\cdots s_{j_{k}}$.
\item
Let $\nset{\infty}
\equiv \{(i_{n})_{n\in {\bf N}}:\mbox{for any } n,\,i_{n}\in\nset{}\}$.
For $J=(j_{n})_{n\in {\bf N}}\in\nset{\infty}$,
$P_{N}(J)$ is a class of $({\cal H},\pi)$ 
with a unit cyclic vector $\Omega\in {\cal H}$ such that
$\{\pi(s_{J_{(n)}})^{*}\Omega:n\in {\bf N}\}$
is an orthonormal family in ${\cal H}$ where
$J_{(n)}\equiv (j_{1},\ldots,j_{n})$. 
\end{enumerate}
The vector $\Omega$ in both (ii) and (iii) is 
called the GP vector of $({\cal H},\pi)$.
\end{defi}
A representation 
$({\cal H},\pi)$ is a {\it {\it cycle} ({\it chain})} if
there is $J\in\nset{k}$ ({\it resp.} $J\in\nset{\infty}$)
such that $({\cal H},\pi)$ is $P_{N}(J)$.
We review results of permutative representations \cite{BJ,DaPi2,DaPi3,K1}.
Any permutative representation is uniquely decomposed into
cyclic permutative representations up to unitary equivalence.
Any cyclic permutative representation is equivalent to $P_{N}(J)$ for a certain 
$J\in \nset{\#}\equiv\coprod_{k\geq 1}\nset{k}\sqcup \nset{\infty}$.
For any $J$, $P_{N}(J)$ exists uniquely up to unitary equivalence.
%
%
\begin{Thm}
\label{Thm:general}
For the operation $\ptimes$ in (\ref{eqn:definition}), the following holds.
\begin{enumerate}
\item
If both $\pi_{1}$ and $\pi_{2}$ are permutative representations, 
then $\pi_{1}\ptimes \pi_{2}$ is also a permutative representation.
\item
If both $\pi_{1}$ and $\pi_{2}$ are cycles, 
then $\pi_{1}\ptimes \pi_{2}$ is a direct sum of cycles.
\item
If $\pi_{1}$ is a permutative representation and $\pi_{2}$ is a chain, 
then $\pi_{1}\ptimes \pi_{2}$ is a direct sum of chains.
\end{enumerate}
\end{Thm}

\noindent
By Lemma \ref{lem:fundamental} (i) and Theorem \ref{Thm:general} (i),
the tensor product of permutative representations is decomposed into cyclic
permutative representations uniquely up to unitary equivalence.
Hence, the decomposition formula of tensor products makes sense.

In order to show more detail, we prepare several notions of multiindices.
Define $\nset{*}_{1}\equiv\coprod_{k\geq 1}\nset{k}$ and
$\nset{*}\equiv\coprod_{k\geq 0}\nset{k}$, $\nset{0}\equiv \{0\}$.
The {\it length} $|J|$ of $J\in \nset{\#}$ is defined by
$|J|\equiv k$ when $J\in \nset{k}$.
For $J_{1},J_{2}\in\nset{*}$ and $J_{3}\in\nset{\infty}$,
$J_{1}\cup J_{2}\equiv(j_{1},\ldots,j_{k},j_{1}^{'},\ldots,j_{l}^{'})$,
$J_{1}\cup J_{3}\equiv(j_{1},\ldots,j_{k},j_{1}^{''},j_{2}^{''},\ldots)$
for $J_{1}=(j_{a})_{a=1}^{k}$, $J_{2}=(j_{b}^{'})_{b=1}^{l}$
and $J_{3}=(j_{n}^{''})_{n\in {\bf N}}$.
Especially, we define $J\cup (0)=(0)\cup J=J$ for convention.
For $J\in\nset{*}$ and $k\geq 2$,
$J^{k}\equiv J\cup\cdots\cup J$ ($k$ times).

For $a,b\in {\bf N}$, we denote 
the greatest common divisor and the least common multiple
of $a$ and $b$ by ${\rm gcd}(a,b)$ and ${\rm lcm}(a,b)$, respectively.
We generalize ${\rm lcm}(a,b)$ as
${\rm lcm}(a,b)=\infty$ when $a=\infty$ or $b=\infty$.
In order to describe the decomposition formula,
we introduce two products among multiindices.
For $K=(k_{i})_{i=1}^{a}\in\nset{a}$ and 
$L=(l_{i})_{i=1}^{a}\in\{1,\ldots,M\}^{a}$ with $1\leq a\leq \infty$,
define $K\cdot L=(x_{i})_{i=1}^{a}\in \{1,\ldots,NM\}^{a}$ by
\[x_{i}\equiv M(k_{i}-1)+l_{i}\quad(i=1,\ldots,a).\]
For $K\in\nset{a}$, $L\in\{1,\ldots,M\}^{b}$ and $C\equiv {\rm lcm}(a,b)$,
define $K*L\in \{1,\ldots,NM\}^{C}$ by
\[K*L\equiv 
\left\{
\begin{array}{ll}
K^{m_{1}}\cdot L^{m_{2}}\quad &(a,b<\infty,\,
m_{1}\equiv C/a,\, m_{2}\equiv C/b),\\
\\
K\cdot L^{\infty}\quad &(a=\infty,\, b<\infty),\\
\\
K\cdot L\quad &(a=b=\infty).\\
\end{array}
\right.
\]
When representations $\pi_{1}$ and $\pi_{2}$ are $P_{N}(K)$ and 
$P_{M}(L)$, we denote $\pi_{1}\ptimes \pi_{2}$ by
$P_{N}(K)\ptimes P_{M}(L)$ simply.
%
%
\begin{Thm}
\label{Thm:main}
If $K\in\nset{a}$ and $L=(l_{j})_{j=1}^{b}\in\{1,\ldots,M\}^{b}$, then
\[
P_{N}(K)\ptimes P_{M}(L)=
\left\{
\begin{array}{ll}
\disp{
\bigoplus_{i=1}^{{\rm gcd}(a,b)}P_{NM}(K*L^{(i)})\quad}&(a,b<\infty),\\
\\
\disp{\bigoplus_{i=1}^{b}P_{NM}(K*L^{(i)})\quad}&(a=\infty,\,b<\infty),\\
\\
\disp{\bigoplus_{i\in {\bf Z}}P_{NM}(K*L^{(i)})\quad}&(a=b=\infty)\\
\end{array}
\right.
\]
where $L^{(i)}\equiv (l_{i},\ldots,l_{b},l_{1},\ldots,l_{i-1})$
when $b<\infty$ and $L^{(i)}\equiv (l_{j}^{'})_{j=1}^{\infty}$
for $l_{j}^{'}\equiv l_{j+i}$ when $j+i\geq 1$ and
$l_{j}^{'}\equiv 1$ when $j+i\leq 0$.
\end{Thm}

In $\S$ \ref{section:second}, we prove Theorem \ref{Thm:symmetry},
Theorem \ref{Thm:general} and Theorem \ref{Thm:main}.
In $\S$ \ref{section:third},
we define a tensor product among certain endomorphisms and states.
In $\S$ \ref{section:fourth}, we show examples and an application
of Theorem \ref{Thm:main}.
In $\S$ \ref{subsection:fourthtwo},
we give two inequivalent endomorphisms of $\co{4}$ and 
show that they are irreducible and do not have inverse.

%
%
\sftt{Proof of Theorems}
\label{section:second}
For $N\geq 2$, let $\con$ be the {\it Cuntz algebra} \cite{C}, that is,
a C$^{*}$-algebra which is universally generated by
generators $s_{1},\ldots,s_{N}$ satisfying
$s_{i}^{*}s_{j}=\delta_{ij}I$ for $i,j\edot$ and
$s_{1}s_{1}^{*}+\cdots +s_{N}s_{N}^{*}=I$.
\\
%
%

\noindent
{\it Proof of Theorem \ref{Thm:symmetry}.}
Let $t_{1}^{(i)},\ldots,t_{N}^{(i)}$ be canonical generators of $\co{N_{i}}$
and let $({\cal H}_{i},\pi_{i})\in {\rm Rep}\co{N_{i}}$ for $i=1,\ldots,n$.
For $\sigma\in {\goth S}_{n}$,
define $\Pi_{\sigma}\equiv \pi_{\sigma(1)}\ptimes
\cdots \ptimes \pi_{\sigma(n)}$,
${\cal K}_{\sigma}\equiv {\cal H}_{\sigma(1)}\otimes
\cdots\otimes {\cal H}_{\sigma(n)}$,
$D_{\sigma}\equiv \co{N_{\sigma(1)}}\otimes
\cdots \otimes  \co{N_{\sigma(n)}}$ and
$\Phi_{\sigma}(s_{[i_{1},\ldots,i_{n}]_{\sigma}})
\equiv t_{i_{\sigma(1)}}^{(\sigma(1))}\otimes
\cdots \otimes t_{i_{\sigma(n)}}^{(\sigma(n))}$. 
For $\sigma,\eta\in {\goth S}_{n}$,
define the unitary $U_{\sigma,\eta}$ from 
${\cal K}_{\eta}$ to ${\cal K}_{\sigma}$ by
$U_{\sigma,\eta}(v_{\eta(1)}\otimes \cdots\otimes v_{\eta(n)})\equiv 
v_{\sigma(1)}\otimes \cdots\otimes v_{\sigma(n)}$
for $v_{1}\otimes \cdots\otimes v_{n}\in {\cal K}_{id}$
and define the map $\tau_{\sigma,\eta}$ from
$D_{\eta}$ to $D_{\sigma}$ by
$\tau_{\sigma,\eta}(x_{\eta(1)}\otimes \cdots \otimes x_{\eta(n)})
\equiv x_{\sigma(1)}\otimes \cdots \otimes x_{\sigma(n)}$
for $x_{1}\otimes \cdots \otimes x_{n}\in D_{id}$.
Then $\pi_{\sigma(1)}\otimes\cdots \otimes \pi_{\sigma(n)}
={\rm Ad}U_{\sigma,\eta}\circ (\pi_{\eta(1)}\otimes
\cdots \otimes \pi_{\eta(n)})\circ \tau_{\eta,\sigma}$
and
$\tau_{\eta,\sigma}\circ \Phi_{\sigma}
=\Phi_{\eta}\circ \alpha_{\eta,\sigma}$.
Then we see that
$\Pi_{\sigma}
=(\pi_{\sigma(1)}\otimes\cdots \otimes \pi_{\sigma(n)})\circ \Phi_{\sigma}
= {\rm Ad}U_{\sigma,\eta}\circ \Pi_{\eta}\circ \alpha_{\eta,\sigma}$.
Hence the statement holds.
\qedh

\noindent
Especially,  the following holds.
%
%
\begin{cor}
\label{cor:transposition}
For $({\cal H}_{1},\pi_{1})\in {\rm Rep}\con$ and $({\cal H}_{2},\pi_{2})
\in {\rm Rep}\co{M}$,
define $\alpha\in {\rm Aut}\co{NM}$ by $\alpha(s_{N(j-1)+i})\equiv s_{M(i-1)+j}$
for $i=1,\ldots,N$ and $j=1,\ldots,M$,
and define the unitary $U$ from ${\cal H}_{1}\otimes {\cal H}_{2}$ to
${\cal H}_{2}\otimes {\cal H}_{1}$ by $U(x\otimes y)\equiv y\otimes x$
for $x\in {\cal H}_{1}$ and $y\in {\cal H}_{2}$.
Then the following holds:
\[\pi_{2}\ptimes \pi_{1}={\rm Ad}U\circ (\pi_{1}\ptimes \pi_{2})\circ \alpha.\]
\end{cor}

\noindent
In this sense, the transposition of tensor product
is understood by the permutation of canonical generators of $\co{NM}$
up to unitary equivalence.
In consequence, $\pi_{1}\ptimes \pi_{2}\sim \pi_{2}\ptimes \pi_{1}$ on 
the fixed-point subalgebra $(\co{NM})^{\alpha}$.

Let $({\cal H},\pi)$ be $P_{N}(J)$ with the GP vector $\Omega$
for $J=(j_{1},\ldots,j_{k})$.
Define $\Omega_{1},\ldots,\Omega_{k}\in {\cal H}$
by
$\Omega_{i}\equiv \pi(s_{j_{i}}\cdots s_{j_{k}})\Omega$
for $i=1,\ldots,k$.
We call $\{\Omega_{i}\}_{i=1}^{k}$ the cycle of vectors of 
$({\cal H},\pi)$.
Let $({\cal H},\pi)$ be $P_{N}(J)$ with the GP vector $\Omega$
for $J=(j_{n})_{n\in {\bf N}}$.
Define the family $\{\Omega_{n}\}_{n\in {\bf N}}$
of vectors in ${\cal H}$ by
$\Omega_{n}\equiv \pi(s_{j_{1}}\cdots s_{j_{n}})^{*}\Omega$
for $n\in {\bf N}$.
We call $\{\Omega_{n}\}_{n\in {\bf N}}$ the chain of vectors of 
$({\cal H},\pi)$.
For a subset $S$ of a Hilbert space ${\cal H}$,
we denote $S^{\perp}$ the orthogonal compliment of $S$ in ${\cal H}$.
\\

\noindent
{\it Proof of Theorem \ref{Thm:general}.}
(i) By definition of the permutative representation
and the tensor product in (\ref{eqn:standard}), statement holds.

\noindent
(ii)
Let $({\cal H}_{1},\pi_{1})$ and $({\cal H}_{2},\pi_{2})$ 
be $P_{N}(J)$ and $P_{M}(K)$, 
and let $\Omega_{1},\ldots \Omega_{a}$ and
$\Omega_{1}^{'},\ldots \Omega_{b}^{'}$ be their cycles of vectors.
Define $\Omega_{i,j}\equiv \Omega_{i}\otimes \Omega_{j}^{'}$ for 
$i=1,\ldots,a$ and $j=1,\ldots,b$.
Then for any $i,j$,
there are $c\geq 1$ and $L\in\{1,\ldots,M\}^{c}$ such that
$(\pi_{1}\ptimes \pi_{2})(s_{L})\Omega_{i,j}=\Omega_{i,j}$.
On the other hand, there exist neither cycle nor chain in both
$\{\Omega_{1},\ldots,\Omega_{a}\}^{\perp}$ and
$\{\Omega_{1}^{'},\ldots,\Omega_{b}^{'}\}^{\perp}$
with respect to $\pi_{1}$ and $\pi_{2}$, respectively.
This implies that there exist neither cycle nor chain in
$\{\Omega_{i,j}:i=1,\ldots,a,\,j=1,\ldots,b\}^{\perp}$.
In consequence, the statement holds.

\noindent
(iii) (a) Let $J\in\nset{a}$ and $K\in\{1,\ldots,M\}^{\infty}$.
Let $({\cal H}_{1},\pi_{1})$  and
$({\cal H}_{2},\pi_{2})$ be $P_{N}(J)$ and $P_{M}(K)$ and
let $\Omega_{1},\ldots \Omega_{a}$ be cycles of vectors of $\pi_{1}$ and
$\Omega_{n}^{'}$, $n\in {\bf N}$ be chain of vectors of $\pi_{2}$.
Define $\Omega_{i,j}\equiv \Omega_{i}\otimes \Omega_{j}^{'}$ for 
$i=1,\ldots,a$ and $j\in {\bf N}$.
Then for any $i,j$, 
there are $i^{'},j^{'}$ and $L\in\{1,\ldots,M\}^{c}$ such that
$(\pi_{1}\ptimes \pi_{2})(s_{L}^{*})\Omega_{i,j}\ne \Omega_{i^{'},j^{'}}$.
On the other hand,
there exists neither cycle nor chain in 
$\{\Omega_{1},\ldots,\Omega_{a}\}^{\perp}$ and
${\cal H}_{2}$ has a chain 
$\{\Omega_{n}^{'}:n\in {\bf N}\}$ and does not have cycle.
This implies that there is neither cycle nor chain in
$(\{\Omega_{i}:i=1,\ldots,a\}\otimes {\cal H}_{2})^{\perp}$.
In consequence, the statement holds.

\noindent
(b) For $J,K\in\{1,\ldots,M\}^{\infty}$,
the statement holds by the similarity of (a).
\qedh

\noindent
{\it Proof of Theorem \ref{Thm:main}.}
We identify $\co{NM}$ as a subalgebra $\varphi_{N,M}(\co{NM})$
of $\con\otimes \co{M}$.
Let $c\equiv {\rm gcd}(a,b)$, $C\equiv {\rm lcm}(a,b)$,
$m_{1}\equiv C/a$ and $m_{2}\equiv C/b$.
Let $({\cal H}_{1},\pi_{1})$ and $({\cal H}_{2},\pi_{2})$ be 
$P_{N}(K)$ and $P_{M}(L)$ with GP vectors $\Omega_{1}$ and $\Omega_{2}$, 
respectively.
Define $\Pi\equiv \pi_{1}\ptimes\pi_{2}$ and
${\cal K}\equiv {\cal H}_{1}\otimes {\cal H}_{2}$.

Assume that $a,b<\infty$. Define 
\[u_{i}\equiv \pi_{1}(t_{k_{i}}\cdots t_{k_{a}})\Omega_{1},\quad
v_{j}\equiv \pi_{2}(r_{l_{j}}\cdots r_{l_{b}})\Omega_{2},\quad
w_{i,j} \equiv u_{i}\otimes v_{j}\]
for $i=1,\ldots,a$ and $j=1,\ldots,b$.
Then
\[\pi_{1}(t_{K})u_{1}=u_{1},\quad
\pi_{2}(r_{L^{(i)}})v_{i}=v_{i}\quad(i=1,\ldots,b).\]
By assumption,
$W\equiv \{w_{i,j}:(i,j)\in\{1,\ldots,a\}\times \{1,\ldots,b\}\}$
is an orthonormal family in ${\cal K}$ and
there exists neither cycle nor chain of vectors in $W^{\perp}$.
Therefore it is sufficient to check cycles in $W$ with respect to $\Pi$
by Theorem \ref{Thm:general} (ii).
Because $\Pi(s_{M(i-1)+j})=\pi_{1}(t_{i})\otimes \pi_{2}(r_{j})$,
$\Pi(s_{K*L^{(i)}})=\pi_{1}(t_{K})^{m_{1}}\otimes\pi_{2}(r_{L^{(i)}})^{m_{2}}$.
These imply that 
\[\Pi(s_{K*L^{(i)}})w_{1,i}=w_{1,i}\quad (i=1,\ldots,c).\]
Define $V_{i}\equiv \overline{\Pi(\co{NM})w_{1,i}}$ for $i=1,\ldots,c$.
For $x,y\in {\bf N}$, 
let $w_{x,y}=w_{i,j}$ if $x\equiv i$  mod $a$ and $y\equiv j$  mod $b$.
Then $w_{1,i},w_{2,i+1},\ldots,w_{1+C-1,i+C-1}$ belong to $V_{i}$
and they are orthogonal.
From this, $W\subset V_{1}\oplus \cdots \oplus V_{c}$.
This implies that ${\cal K}=V_{1}\oplus \cdots \oplus V_{c}$
and $(V_{i},\Pi|_{V_{i}})$ is $P_{NM}(K*L^{(i)})$.
Hence the statement holds for this case.

Assume that $a=\infty$ and $b<\infty$.
Define 
\[u_{1}\equiv \Omega_{1},\quad
u_{i}\equiv \pi_{1}(t_{k_{1}}\cdots t_{k_{i-1}})^{*}\Omega_{1},\quad
v_{j}\equiv \pi_{2}(r_{l_{j}}\cdots r_{l_{b}})\Omega_{2},\quad
w_{i,j} \equiv u_{i}\otimes v_{j}\]
for $i\geq 1$ and $j=1,\ldots,b$.
Then $\pi_{2}(r_{L^{(i)}})v_{i}=v_{i}$ for $i=1,\ldots,b$.
By assumption,
$W\equiv \{w_{i,j}:i\geq 1,\, j\in \{1,\ldots,b\}\}$
is an orthonormal family in ${\cal K}$ and
there exists neither cycle nor chain in $W^{\perp}$.
Therefore it is sufficient to check chains in $W$ with respect to $\Pi$
by Theorem \ref{Thm:general} (iii).
Denote $K*L^{(i)}=(d_{i,j})_{j=1}^{\infty}$ and $(L^{(i)})^{\infty}
=(l_{i,j}^{'})_{j=1}^{\infty}$.
Then $d_{i,j}=M(k_{j}-1)+l_{i,j}^{'}$.
These imply that 
\[\Pi(s_{d_{i,1}}\cdots s_{d_{i,mb}})^{*}w_{1,i}=w_{1+mb,i}
\quad (i=1,\ldots,b,\,m\geq 1).\]
Define $V_{i}\equiv \overline{\Pi(\co{NM})w_{1,i}}$ for $i=1,\ldots,b$.
For $x,y\in {\bf N}$,  let $w_{x,y}=w_{x,j}$ if $y\equiv j$ mod $b$.
Then $w_{n,n+i-1}$ belongs to $V_{i}$ for $n\geq 1$ and they are orthogonal
and $w_{n,n+i-1+bm}=w_{n,n+i-1}$ for each $n,m\geq 1$ and $i=1,\ldots b$.
From this, $W\subset V_{1}\oplus \cdots \oplus V_{b}$.
This implies that ${\cal K}=V_{1}\oplus \cdots \oplus V_{b}$
and $(V_{i},\Pi|_{V_{i}})$ is $P_{NM}(K*L^{(i)})$.
Hence the statement holds for this case.

Assume that $a=b=\infty$. Define 
\[u_{1}\equiv \Omega_{1},\quad
u_{n}\equiv \pi_{1}(t_{k_{1}}\cdots t_{k_{n-1}})^{*}\Omega_{1}\quad(n\geq 2),\]
\[v_{m}\equiv \pi_{2}(r_{l_{1}}\cdots r_{l_{m-1}})^{*}\Omega_{2},\quad
v_{1}\equiv \Omega_{2},\quad
v_{-m+2}\equiv \pi_{2}(r_{1}^{m-1})\Omega_{2},\quad(m\geq 2),\]
\[w_{n,m} \equiv u_{n}\otimes v_{m}\quad ((n,m)\in{\bf N}\times {\bf Z}).\]
By assumption, $W\equiv \{w_{n,m}:(n,m)\in{\bf N}\times {\bf Z}\}$ 
is an orthonormal family in ${\cal K}$ 
and there exists neither cycle nor chain in $W^{\perp}$.
Therefore it is sufficient to check chains in $W$ with respect to $\Pi$
by Theorem \ref{Thm:general} (iii).
Denote $K*L^{(i)}=(d_{i,j})_{j=1}^{\infty}$ and $(L^{(i)})^{\infty}
=(l_{i,j}^{'})_{j=1}^{\infty}$.
Then $d_{i,j}=M(k_{j}-1)+l_{i,j}^{'}$.
These imply that 
\[\Pi(s_{d_{i,1}}\cdots s_{d_{i,m}})^{*}w_{1,1+i}
=w_{1+m,1+i+m}\quad ((i,m)\in{\bf Z}\times {\bf N}).\]
Define $V_{i}\equiv \overline{\Pi(\co{NM})w_{1,i}}$ for $i\in {\bf Z}$.
Then $w_{n,n+i}$ belongs to $V_{i}$ for $n\geq 1$ and they are orthogonal.
From this, $W\subset \bigoplus_{i\in {\bf Z}}V_{i}$.
This implies that ${\cal K}=\bigoplus_{i\in {\bf Z}}V_{i}$
and $(V_{i},\Pi|_{V_{i}})$ is $P_{NM}(K*L^{(i)})$.
Hence the statement holds for this case.
\qedh

%
%
\sftt{Tensor product of others}
\label{section:third}

We introduce a tensor product among
certain endomorphisms and among states
such that it is compatible with the $\varphi$-tensor product
in (\ref{eqn:definition}).

%
%
\ssft{Tensor product of endomorphisms}
\label{subsection:thirdone}
Let $V_{N,l}\equiv ({\bf C}^{N})^{\otimes l}$
and let $\{\vep_{i}\}_{i=1}^{N}$ be the standard basis of ${\bf C}^{N}$.
Define $\vep_{J}\equiv \vep_{j_{1}}\otimes \cdots\otimes \vep_{j_{l}}$
for $J=(j_{1},\ldots,j_{l})\in\nset{l}$.
Define $U(N,l)$ the group of all unitaries on $V_{N,l}$
with respect to the standard inner product $\langle\cdot|\cdot\rangle$
of $V_{N,l}$.
For $g\in U(N,l)$, 
define $g_{JK}\equiv \langle \vep_{J}|\,g\,\vep_{K}\rangle $ 
for $J,K\in\nset{l}$ and define the endomorphism $\psi_{g}$ of $\con$ by
%
%
\begin{equation}
\label{eqn:defeqn}
\psi_{g}(s_{i})\equiv u_{g}s_{i}\quad(i\edot)
\end{equation}
where $u_{g}\equiv\sum_{J,K\in \nset{l}}g_{JK}s_{J}(s_{K})^{*}$.
This $\psi_{g}$ is called the {\it generalized
permutative endomorphism} of $\con$ by $g$ \cite{PE01,PE02}.
Especially, if $g\in U(N,1)\cong U(N)$, then 
$\psi_{g}$ is the canonical $U(N)$-action on $\con$.

Let denote ${\rm End}\con$ the set of all endomorphisms of $\con$.
For $\rho_{1}\in {\rm End}\con$, $\rho_{2}\in {\rm End}\co{M}$ and
$\varphi_{N,M}$ in (\ref{eqn:standard}),
consider $(\rho_{1}\otimes \rho_{2})\circ \varphi_{N,M}\in 
{\rm Hom}(\co{NM},\con\otimes \co{M})$.
Remark that
$\{(\rho_{1}\otimes \rho_{2})\circ \varphi_{N,M}\}(\co{NM})$
is not a subset of $\varphi_{N,M}(\co{NM})$ in general.
%
%
\begin{prop}
\label{prop:morphism}
\begin{enumerate}
\item
If $g\in U(N,l)$ and $h\in U(M,k)$, then 
$\{(\psi_{g}\otimes \psi_{h})\circ \varphi_{N,M}\}(\co{NM})\subset
\varphi_{N,M}(\co{NM})$.
\item
For $g\in U(N,l)$ and $h\in U(M,k)$, 
define $\psi_{g}\ptimes \psi_{h}\in {\rm End}\co{NM}$ by
%
%
\begin{equation}
\label{eqn:morphism}
\psi_{g}\ptimes \psi_{h}\equiv 
(\varphi_{N,M})^{-1}\circ (\psi_{g}\otimes \psi_{h})\circ \varphi_{N,M}.
\end{equation}
Then there exists $\hat{g}\in U(NM,m)$ for $m\equiv \max\{l,k\}$ such that
$\psi_{g}\ptimes \psi_{h}=\psi_{\hat{g}}$.
\item
Let $\alpha_{N}$ be the canonical $U(N)$-action on $\con$.
For $g=(g_{ij})\in U(N)$ and $h=(h_{ij})\in U(M)$,
\[\alpha_{N,g}\ptimes \alpha_{M,h}=\alpha_{NM,g*h}\]
where $g*h=((g*h)_{a,b})_{a,b=1}^{NM}\in U(NM)$ is defined by
\[(g*h)_{M(k-1)+l,M(i-1)+j}\equiv g_{ki}h_{lj}
\quad(i,k=1,\ldots,N,\,j,l=1,\ldots,M).\]
\item
For $l\geq 1$, define 
\[UE_{N,l}\equiv \{\psi_{g}:g\in U(N,l)\}.\] 
For $l_{1},l_{2},l_{3}\in {\bf N}$,
if $\rho_{1}\in UE_{N,l_{1}}$, $\rho_{2}\in UE_{M,l_{2}}$
and $\rho_{3}\in UE_{L,l_{3}}$,
then $\rho_{1}\ptimes(\rho_{2}\ptimes\rho_{3})
=(\rho_{1}\ptimes\rho_{2})\ptimes\rho_{3}$.
\item
For any $\pi_{1}\in {\rm Rep}\con$, $\pi_{2}\in {\rm Rep}\co{M}$,
$\rho_{1}\in UE_{N,l}$ and $\rho_{2}\in UE_{M,l}$,
$(\pi_{1}\ptimes \pi_{2})\circ (\rho_{1}\ptimes \rho_{2})
=(\pi_{1}\circ \rho_{1})\ptimes (\pi_{2}\circ \rho_{2})$.
\end{enumerate}
\end{prop}
%
%
\pr
Let $t_{1},\ldots,t_{N},r_{1},\ldots,r_{M}$ and $s_{1},\ldots,s_{NM}$
be canonical generators of $\con$, $\co{M}$ and $\co{NM}$, respectively.
When $l=k$, we see that $u_{g}\otimes u_{h}\in \varphi_{N,M}(\co{NM})$.
Define $U\equiv \varphi^{-1}_{N,M}(u_{g}\otimes u_{h})$. Then 
$(\psi_{g}\ptimes \psi_{h})(s_{i})=Us_{i}$ for $i=1,\ldots,NM$.
Assume that $m=l>  k$.
Replace $h$ by $\hat{h}\in U(M,l)$ defined by $\hat{h}\equiv h\otimes I$.
Then $u_{\hat{h}}=u_{h}$ and $\psi_{\hat{h}}=\psi_{h}$.
Hence $u_{g}\otimes u_{h}=u_{g}\otimes u_{\hat{h}}
\in \varphi_{N,M}(\co{NM})$.
Define $U\equiv \varphi_{N,M}^{-1}(u_{g}\otimes u_{\hat{h}})$.
Then
\[\varphi_{N,M}(U)=u_{g}\otimes u_{\hat{h}}
=\sum_{J,K\in\nset{l}}\sum_{P,Q\in\{1,\ldots,M\}^{l}}
g_{JK}\hat{h}_{PQ}\cdot t_{J}t_{K}^{*}\otimes r_{P}r_{Q}^{*}.\]
Define 
\[\phi_{N,M}(i,j)\equiv M(i-1)+j,\quad
\phi_{N,M}(J,K)\equiv(\phi_{N,M}(j_{1},k_{1}),\ldots,\phi_{N,M}(j_{m},k_{m}))\]
for $J=(j_{1},\ldots,j_{m})\in\nset{m}$ and 
$K=(k_{1},\ldots,k_{m})\in\{1,\ldots,M\}^{m}$.
Because $\varphi_{N,M}(s_{\phi_{N,M}(i,j)})=t_{i}\otimes r_{j}$,
\[t_{J}t_{K}^{*}\otimes r_{P}r_{Q}^{*}=
\varphi_{N,M}(s_{\phi_{N,M}(J,P)}s_{\phi_{N,M}(K,Q)}^{*}).\]
From this, we see that
%
%
\begin{equation}
\label{eqn:endomorphism}
(\psi_{g}\ptimes \psi_{h})(s_{M(i-1)+j})=Us_{M(i-1)+j}
\quad(i=1,\ldots,N,\,j=1,\ldots,M).
\end{equation}
By definition of $U$, (i) is shown.

On the other hand,
define $g\times_{\phi}h=((g\times_{\phi}h)_{JK})\in U(NM,m)$ by
\[(g\times_{\phi} h)_{\phi_{N,M}(J,P),
\phi_{N,M}(K,Q)}\equiv g_{JK}\hat{h}_{PQ}\quad (J,K\in\{1,\ldots,NM\}^{m}).\]
Then we can verify that $U=u_{g\times_{\phi}h}$.
By (\ref{eqn:endomorphism}), $U$ is a unitary in $\co{NM}$ and
$\psi_{g}\ptimes \psi_{h}=\psi_{g\times_{\phi}\hat{h}}$.
Hence (ii) is proved.
As the similarity, it is shown when $k<l$.
The statement of (iii) is the case $(l,k)=(1,1)$.
(iv) is verified by the weak co-associativity of $\varphi$.
(v) holds by definition.
\qedh

\noindent
By Proposition \ref{prop:morphism} (ii),
if $\rho$ and $\rho^{'}$ are generalized
permutative endomorphisms of $\con$ and $\co{M}$, respectively,
then $\rho\ptimes \rho^{'}$ is also a generalized 
permutative endomorphism of $\co{NM}$.
Furthermore $UE_{N,l}\ptimes UE_{M,l}\subset UE_{NM,l}$
for each $l\geq 1$.

%
%
\ssft{Tensor product of states}
\label{subsection:thirdtwo}
Let denote $(\con)^{*}$ the set of all bounded linear maps on $\con$.
For $(f,g)\in (\con)^{*}\times (\co{M})^{*}$, 
define $f\ptimes g\in (\co{NM})^{*}$ by
\[f\ptimes g\equiv (f\otimes g)\circ \varphi_{N,M}.\]
By the weak co-associativity of $\{\varphi_{N,M}\}_{N,M\geq 2}$,
$(f\ptimes g)\ptimes h =f\ptimes (g\ptimes h)$
for each $f\in(\con)^{*},
g\in(\co{M})^{*}$ and $h\in(\co{L})^{*}$.
Especially,
if $f$ and $g$ are states of $\con$ and $\co{M}$, respectively,
then $f\ptimes g$ is also a state of $\co{NM}$
because $\varphi_{N,M}$ is unital and $*$-preserving.
Let$({\cal H}_{1},\pi_{1})$ and $({\cal H}_{2},\pi_{2})$
be representations of $\con$ and $\co{M}$, respectively.
For vector states $\omega_{1}$ and 
$\omega_{2}$ of $\con$ and $\co{M}$ defined by
$\omega_{1}=\langle \Omega_{1}|\pi_{1}(\cdot)\Omega_{1}\rangle $ and
$\omega_{2}=\langle \Omega_{2}|\pi_{2}(\cdot)\Omega_{2}\rangle $,
we see that
\[\langle \Omega_{1}\otimes \Omega_{2}|
(\pi_{1}\ptimes \pi_{2})(\cdot)(\Omega_{1}\otimes \Omega_{2})\rangle 
=\omega_{1}\otimes_{\varphi}\omega_{2}.\]
In consequence,
the $\varphi$-tensor of two vector states is also a vector state.

%
%
\sftt{Example and application}
\label{section:fourth}
We show examples and an application of Theorem \ref{Thm:main}.
%
%
\ssft{Example}
\label{subsection:fourthone}
For any $i\in\nset{}$ and $j\in\{1,\ldots,M\}$, 
\[P_{N}(i)\ptimes P_{M}(j)=P_{NM}(M(i-1)+j).\]
For example, when $N=M=2$,
\[P_{2}(1)\ptimes P_{2}(1)=P_{4}(1),\quad P_{2}(1)\ptimes P_{2}(2)=P_{4}(2),\]
\[P_{2}(2)\ptimes P_{2}(1)=P_{4}(3),\quad P_{2}(2)\ptimes P_{2}(2)=P_{4}(4).\]
By Theorem 2.7 (iv) in \cite{PE01},
$P_{2}(1)\ptimes P_{2}(2)=P_{4}(2)\not \sim P_{4}(3)=P_{2}(2)\ptimes P_{2}(1)$.
Hence $\ptimes$ among representations is not symmetric in general.

For $i,j\in\nset{}$ and $i^{'},j^{'}\in\{1,\ldots,M\}$, 
\[
\begin{array}{rl}
P_{N}(ij)\ptimes P_{M}(i^{'}j^{'})=&
P_{NM}(M(i-1)+i^{'},M(j-1)+j^{'})\\
&\oplus P_{NM}(M(i-1)+j^{'},M(j-1)+i^{'}).\\
\end{array}
\]

For a representation $\pi$,
we denote $\pi^{\ptimes n}$ the $n$-times $\varphi$-tensor product of $\pi$.
By the induction method, we have the following.
\[P_{2}(12)^{\ptimes n}=\bigoplus_{i=1}^{2^{n-1}}P_{2^{n}}(i,2^{n}-i+1)
\quad(n\geq 1).\]
Especially, this decomposition is the irreducible decomposition 
and multiplicity-free.
For example,
\[P_{2}(12)^{\ptimes 3}
=P_{8}(18)\oplus P_{8}(36)\oplus P_{8}(27)\oplus P_{8}(45).\]

%
%
\ssft{Application ---Two endomorphisms of $\co{4}$}
\label{subsection:fourthtwo}
Let ${\cal A}$ be a unital C$^{*}$-algebra ${\cal A}$.
An endomorphism $\rho$ of ${\cal A}$ is {\it irreducible} if 
$\rho({\cal A})^{'}\cap {\cal A}={\bf C}I$.
Two endomorphisms  $\rho$ and $\rho^{'}$ are {\it equivalent}
if there exists a unitary $u\in {\cal A}$ such that 
${\rm Ad}u\circ \rho=\rho^{'}$.
In algebraic quantum field theory \cite{DHR},
equivalence classes of irreducible endomorphisms are important.
We show such examples.
Let $s_{1},s_{2},s_{3},s_{4}$ be canonical generators of $\co{4}$.
%
%
\begin{prop}
\label{prop:branching}
Define $\rho,\bar{\rho}\in {\rm End}\co{4}$ by
\[
\begin{array}{rl}
\rho(s_{1})\equiv&s_{23,1}+s_{21,2}+s_{14,3}+s_{12,4},\\
\rho(s_{2})\equiv&s_{13,1}+s_{11,2}+s_{24,3}+s_{22,4},\\
\rho(s_{3})\equiv&s_{41,1}+s_{43,2}+s_{32,3}+s_{34,4},\\
\rho(s_{4})\equiv&s_{31,1}+s_{33,2}+s_{42,3}+s_{44,4},\\
\end{array}
\quad
\begin{array}{rl}
\bar{\rho}(s_{1})\equiv&s_{32,1}+s_{14,2}+s_{31,3}+s_{13,4},\\
\bar{\rho}(s_{2})\equiv&s_{41,1}+s_{23,2}+s_{42,3}+s_{24,4},\\
\bar{\rho}(s_{3})\equiv&s_{12,1}+s_{34,2}+s_{11,3}+s_{33,4},\\
\bar{\rho}(s_{4})\equiv&s_{21,1}+s_{43,2}+s_{22,3}+s_{44,4}\\
\end{array}
\]
where $s_{ij,k}\equiv s_{i}s_{j}s_{k}^{*}$ for $i,j,k= 1,2,3,4$.
Then the following holds.
\begin{enumerate}
\item
$P_{4}(1)\circ \rho=P_{4}(24)$ and 
$P_{4}(1)\circ \bar{\rho}=P_{4}(34)$.
\item
Both $\rho$ and $\bar{\rho}$ are irreducible
and not automorphisms.
\item
$\rho$ and $\bar{\rho}$ are are not equivalent.
\end{enumerate}
\end{prop}
%
%
\pr
Let $t_{1},t_{2}$ be canonical generators of $\co{2}$.
Define $\psi_{12},\psi_{13}\in {\rm End}\co{2}$
by
\[\psi_{12}(t_{1})\equiv t_{12,1}+t_{11,2},\,\,
\psi_{12}(t_{2})\equiv t_{2},\quad
\psi_{13}(t_{1})\equiv t_{21,1}+t_{12,2},\,\,
\psi_{13}(t_{2})\equiv  t_{11,1}+t_{22,2}.\]
We can verify that
$\rho= \psi_{12}\ptimes\psi_{13}$
and $\bar{\rho}=\psi_{13}\ptimes \psi_{12}$.

\noindent
(i)
By Table II in \cite{PE01}, 
$ P_{2}(1)\circ \psi_{12}=P_{2}(12)$ and 
$P_{2}(1)\circ\psi_{13}=P_{2}(2)$.
From these,
\[
\begin{array}{rll}
P_{4}(1)\circ \rho=&(P_{2}(1)\ptimes P_{2}(1))\circ(
\psi_{12}\ptimes\psi_{13})&(\mbox{by Theorem \ref{Thm:main}})\\
= &(P_{2}(1)\circ\psi_{12})\ptimes (P_{2}(1)\circ\psi_{13})
&(\mbox{by Proposition \ref{prop:morphism} (v)})\\
= &P_{2}(12)\ptimes P_{2}(2)\\
= &P_{4}(24)&(\mbox{by Theorem \ref{Thm:main}}).\\
\end{array}
\]
In the same way,
we obtain $P_{4}(1)\circ \bar{\rho}=P_{4}(34)$.

\noindent
(ii)
Every $P_{4}(1)$, $P_{4}(24)$ and $P_{4}(34)$ 
are irreducible by Theorem 2.7 (iii) in \cite{PE01}.
Hence  both $\rho$ and $\bar{\rho}$ are 
irreducible by Lemma 3.1 (i) in \cite{PE01}.

Define $\kappa_{1},\kappa_{2}\in {\rm Aut}\co{4}$ by
\[(\kappa_{1}(s_{1}),\ldots,\kappa_{1}(s_{4}))\equiv 
(-s_{2},-s_{1},s_{4},s_{3}),\quad
(\kappa_{2}(s_{1}),\ldots,\kappa_{2}(s_{4}))\equiv 
(-s_{3},s_{4},-s_{1},s_{2}).\]
Then $\kappa_{1}\circ \rho=\rho$ and
$\kappa_{2}\circ \bar{\rho}=\bar{\rho}$,
From these,
$\rho(\co{4})\subset \co{4}^{\kappa_{1}}\ne \co{4}$ and
$\bar{\rho}(\co{4})\subset \co{4}^{\kappa_{2}}\ne \co{4}$
where $\co{4}^{\kappa_{i}}
\equiv \{x\in \co{4}:\kappa_{i}(x)=x\}$ for $i=1,2$.
Hence both $\rho$ and $\bar{\rho}$ are not automorphisms.

\noindent
(iii)
By Theorem 2.7 (iv) in \cite{PE01},
$P_{4}(1)\circ \rho=P_{4}(24)\not\sim P_{4}(34)=P_{4}(1)\circ \bar{\rho}$.
Hence $\rho$ and $\bar{\rho}$ are not equivalent
by Lemma 3.1 (ii) in \cite{PE01}.
\qedh

\noindent
By the proof of Proposition \ref{prop:branching},
$\psi_{12}\ptimes \psi_{13}\not \sim 
\psi_{13}\ptimes \psi_{12}$.
Hence $\ptimes$ among endomorphisms is not symmetric in general.
Examples of endomorphisms of $\con$ like $\rho$ and $\bar{\rho}$
in Proposition \ref{prop:branching}
are roughly classified in \cite{PE04} for $N=2,3$.\\

\noindent
{\bf Acknowledgement:}
The author would like to express my sincere thanks to Izumi Ojima 
for his interest in this topic.


\end{document}